\documentclass{tac}

\usepackage{amsmath}

\title{Limit Sketches and the Universal Realization of a
Sketch}

\address{Department of Mathematics, Applied Mathematics, and Statistics\\
Case Western Reserve University\\
 Cleveland, Ohio, 44106}

\author{Johnathon Taylor}

\usepackage{tikz}

\tikzset{
  twocell/.style={double,->,shorten <=2pt,shorten >=2pt}
}

\newcommand{\ob}{\textbf{Ob}}

\newcommand{\Sq}{\textbf{Sq}}

\newcommand{\free}{\mathrm{Fr}}

\newcommand{\limsk}{\mathrm{LS}}

\newcommand{\relimsk}{\mathrm{RLS}}

\begin{document}

\maketitle

\begin{abstract}
We construct the universal realized limit sketch associated to a given limit sketch. The construction uses factorization systems to organize the classical argument of \cite{BastChar}, yielding a streamlined and conceptually unified formulation of the technical steps. This provides a structured framework for understanding realizations of limit sketches in terms of factorization-theoretic data.
\end{abstract}

\copyrightyear{2026}

\keywords{factorization systems, limit sketches}
\amsclass{18A32,18A35,18C30}

\eaddress{jmt240@case.edu}

\section{Introduction}

In Proposition 3 of \cite{BastChar}, Bastiani and Ehresmann construct a universal realized limit sketch associated to a cone-bearing neo-category. This construction specializes to a universal realized limit sketch associated to any limit sketch. The existence of the universal realized limit sketch from a limit sketch is used throughout the literature, such as by Ara in Proposition 1.11 of \cite{Dim1}. Motivated by this, we provide a systematic account of the construction of the universal realized limit sketch.

The construction of the universal realized limit sketch exhibits structural features analogous to a fibrant replacement in a homotopical setting. In this paper, we make this analogy precise in the context of limit sketches. We show that the assignment of a limit sketch to its universal realization extends to a $2$-functor, which is left $2$-adjoint to the forgetful $2$-functor from realized limit sketches to limit sketches. As a consequence, the canonical structure map from a limit sketch to its realization is a Morita equivalence of limit sketches.

The paper is organized as follows. Section 2 recalls background on factorization systems. Sections 3 and 4 discusses limits and the $2$-categories of limit sketches and realized limit sketches. In Section 5, we construct a strong factorization system associated to every limit sketch. Using these factorization systems, we define the universal realization as a left $2$-adjoint to the forgetful $2$-functor. In the final section, we prove that the unit of this $2$-adjunction is a Morita equivalence of limit sketches on every component.

The result of this paper is not unique. The novel aspect of our paper lies in the proof strategy that we implement, i.e., the construction of a strong factorization system for every limit sketch and showing how to construct a left $2$-adjoint from these factorization systems.

\section{Orthogonal and Strong Factorization Systems}

We recall basic facts about orthogonal and strong factorization systems that will be used in the construction of the universal realization of a limit sketch. We assume familiarity with weak factorization systems and functorial factorization systems. 

\begin{definition}\label{uni_fact}
Let $f:x\to y$ and $g:a\to b$ be morphisms in a category $\mathcal{C}$. We write $f \downarrow g$ if every commutative square
\[
\begin{tikzpicture}[node distance=2cm,auto]
\node (A) {$x$};
\node (B)[right of=A] {$a$};
\node (C)[below of=A] {$y$};
\node (D)[right of=C] {$b$};
\draw[->] (A) to node {} (B);
\draw[->] (A) to node[left] {$f$} (C);
\draw[->] (B) to node[right] {$g$} (D);
\draw[->] (C) to node {} (D);
\end{tikzpicture}
\]
admits a unique diagonal filler.
\end{definition}

\begin{notation}
Let $\mathcal{C}$ be a category and $H$ a class of morphisms of $\mathcal{C}$. We define
\[
H^{\downarrow} := \{ f \mid h \downarrow f \text{ for all } h \in H \},
\qquad
H^{\uparrow} := \{ f \mid f \downarrow h \text{ for all } h \in H \}.
\]
\end{notation}

\begin{definition}
A weak factorization system $(\mathcal{L},\mathcal{R})$ on a category $\mathcal{C}$ is an \emph{orthogonal factorization system} if every lifting problem has a unique solution.
\end{definition}

\begin{definition}
A \emph{strong factorization system} on a category $\mathcal{C}$ consists of an orthogonal factorization system $(\mathcal{L},\mathcal{R})$ together with a functorial factorization
\[
F = (L,E,R) : \mathrm{Arr}(\mathcal{C}) \to \mathrm{Arr}(\mathcal{C}) \times_{\mathcal{C}} \mathrm{Arr}(\mathcal{C})
\]
such that $L(f) \in \mathcal{L}$ and $R(f) \in \mathcal{R}$ for all morphisms $f$.
\end{definition}

\begin{theorem}\label{Kelly_OFS}
Let $\mathcal{C}$ be a category with pullbacks and let $\kappa$ be a regular cardinal. Let $H$ be a set of morphisms whose domains and codomains are $\kappa$-compact. Then $(H^{\downarrow\uparrow}, H^{\downarrow})$ is an orthogonal factorization system.
\end{theorem}

\begin{proposition}\label{upgrade to SFS}
Every orthogonal factorization system $(\mathcal{L},\mathcal{R})$ on a category $\mathcal{C}$ induces a functorial factorization
\[
F : \mathrm{Arr}(\mathcal{C}) \to \mathrm{Arr}(\mathcal{C}) \times_{\mathcal{C}} \mathrm{Arr}(\mathcal{C})
\]
with $L(f) \in \mathcal{L}$ and $R(f) \in \mathcal{R}$ for all $f \in \mathcal{C}$.
\end{proposition}

\begin{proof}
Assuming the axiom of choice, choose for each morphism $f$ a factorization
\[
x \xrightarrow{L(f)} E(f) \xrightarrow{R(f)} y.
\]
Functoriality follows from Proposition 2.8 of \cite{Garn3}.
\end{proof}

\begin{corollary}[Kelly's Small Object Argument]\label{Kelly_SOA}
Let $\mathcal{C}$ be locally presentable and let $I$ be a set of morphisms with sequentially small domains and codomains. Then $I$ generates a strong factorization system on $\mathcal{C}$.
\end{corollary}

\begin{proof}
By Theorem \ref{Kelly_OFS}, $(I^{\downarrow\uparrow}, I^{\downarrow})$ is an orthogonal factorization system. The result follows from Proposition \ref{upgrade to SFS}.
\end{proof}

\medskip

We do not use the $k_{+}$-construction of Kelly \cite{Kell} in this exposition. Instead, we rely on the axiom of choice to simplify the presentation. A construction avoiding choice is available in \cite{Kell}.

\section{Diagrams, Cones, and Limits}

We recall basic definitions concerning diagrams and limits.

\begin{definition}
A \emph{diagram} in a category $\mathcal{E}$ is a functor $\phi : I \to \mathcal{E}$ from a small category $I$.
\end{definition}

\begin{definition}
Given an object $x$ of $\mathcal{E}$ and a small category $I$, the \emph{constant diagram} at $x$ is the functor $*_x : I \to \mathcal{E}$ sending every object to $x$ and every morphism to $\mathrm{id}_x$.
\end{definition}

\begin{definition}\label{cone_def}
A \emph{cone} over a diagram $\phi : I \to \mathcal{E}$ consists of an object $x$ of $\mathcal{E}$ together with a natural transformation
\[
\gamma : *_x \Rightarrow \phi.
\]
We refer to $x$ as the vertex of the cone.
\end{definition}

\begin{definition}
Let $C$ and $C'$ be cones over a diagram $\phi : I \to \mathcal{E}$ with vertices $x$ and $x'$. A morphism of cones $C \to C'$ is a morphism $g : x \to x'$ such that
\[
\gamma'_i \circ g = \gamma_i
\quad \text{for all } i \in I.
\]
\end{definition}

\begin{definition}
A cone $C$ over a diagram $\phi : I \to \mathcal{E}$ is a \emph{limit cone} if for every other cone $C'$ over $\phi$ there exists a unique morphism of cones $C' \to C$.
\end{definition}

\begin{example}[Equalizer in the category of groups]
Let $f,g : G \to H$ be group homomorphisms. The \emph{equalizer} of $f$ and $g$ is the subgroup
\[
\mathrm{Eq}(f,g) := \{ x \in G \mid f(x) = g(x) \}.
\]
There is an inclusion homomorphism
\[
e : \mathrm{Eq}(f,g) \hookrightarrow G
\]
given by the inclusion map.

For any group $X$ and homomorphism $h : X \to G$ such that $f \circ h = g \circ h$, there exists a unique homomorphism
\[
\bar{h} : X \to \mathrm{Eq}(f,g)
\]
such that $e \circ \bar{h} = h$.
\end{example}

\section{Limit Sketches}

We fix a categorical framework for specifying limit constraints via distinguished cones.

\begin{definition}\label{def_of_lim_sk}
A \emph{limit sketch} consists of a pair $(E, \mathcal{C})$ where
\begin{itemize}
\item $E$ is a category and
\item $\mathcal{C} = \{C_\alpha\}_{\alpha \in A}$ is a set of specified cones in $E$ indexed by a set $A$.
\end{itemize}
\end{definition}

\begin{definition}
Let $(E, \{C_\alpha\}_{\alpha \in A})$ and $(E', \{C'_\beta\}_{\beta \in B})$ be limit sketches. A \emph{morphism of limit sketches} consists of
\begin{itemize}
\item a functor $F : E \to E'$ and
\item a function $f : A \to B$
\end{itemize}
such that for each $\alpha \in A$, the cone $F(C_\alpha)$ is equal (up to the specified indexing) to $C'_{f(\alpha)}$.
\end{definition}

\begin{definition}
Let $(F,f),(G,g)$ be morphisms of limit sketches. A \emph{2-morphism} $\phi : (F,f) \Rightarrow (G,g)$ is a natural transformation
\[
\phi : F \Rightarrow G
\]
compatible with the indexing data.
\end{definition}

\begin{notation}
We write $\limsk$ to denote the resulting $2$-category of limit sketches. We write $\relimsk \subseteq \limsk$ for the full sub-$2$-category spanned by those sketches in which every specified cone is a limit cone in $E$. Objects of $\relimsk$ are called \emph{realized limit sketches}.
\end{notation}

\medskip

We will usually suppress explicit reference to the indexing set of cones when no ambiguity arises.

\begin{example}\label{sk_of_fin_sets}
Let $\mathbf{Fin}$ denote a skeleton of finite sets. Then $\mathbf{Fin}^{op}$ is the free finite-product category on one generator, hence determines a realized limit sketch. More generally, every Lawvere theory admits the structure of a realized limit sketch.
\end{example}

\begin{example}
The simplex category $\Delta$, equipped with no specified cones, is a limit sketch in which the cone data is empty.
\end{example}

\begin{notation}\label{und_cat}
The \emph{underlying category} functor
\[
U : \limsk \to \mathbf{Cat}
\]
is defined by sending a limit sketch $(E,\mathcal{C})$ to its underlying category $E$.
\end{notation}

\begin{lemma}
The assignment $U$ extends to a $2$-functor
\[
U : \limsk \to \mathbf{Cat}.
\]
\end{lemma}

\begin{notation}\label{free_lim_sk}
The \emph{free limit sketch} on a category $E$ is the limit sketch
\[
\free(E) := (E, \emptyset),
\]
with no specified cones.
\end{notation}

\begin{lemma}
The assignment $E \mapsto \free(E)$ extends to a $2$-functor
\[
\free : \mathbf{Cat} \to \limsk.
\]
\end{lemma}

\begin{lemma}
The functor 
\[
\free : \mathbf{Cat} \to \limsk.
\]
is left $2$-adjoint to the forgetful $2$-functor
\[
U : \limsk \to \mathbf{Cat}.
\]
\end{lemma}

\section{Generating a Strong Factorization System}

Fix a limit sketch
\[
(E,A,\{C_\alpha\}_{\alpha \in A}),
\]
and an object $y$ of $E$. For each $\alpha \in A$, the associated cone
\[
C_\alpha = (x_\alpha, I_\alpha, \phi_\alpha, \delta_\alpha),
\]
consists of an object $x_\alpha$ of $E$, a small indexing category $I_\alpha$, a diagram $\phi_\alpha : I_\alpha \to E$, and
a natural transformation
\[
\delta_\alpha : *_x \Rightarrow \phi_\alpha.
\]
\medskip

The goal of this section is to construct a strong factorization system on $E/\mathbf{Cat}$ generated by the cone data of the sketch.

\begin{definition}
Let $y$ be an object of $E$ and $\alpha \in A$. A \emph{$(y,\alpha)$-extension square} consists of a category $D$, a functor $F : E \to D$, and a natural transformation $\psi$ fitting into a diagram

\[
\begin{tikzpicture}[>=stealth, thick, node distance=3cm]

\node (A) at (0,0) {$I_\alpha$};
\node (B) at (4,0) {$E$};
\node (C) at (0,-3) {$E$};
\node (D) at (4,-3) {$D$};

\draw[->] (A) to node[above] {$*_y$} (B);
\draw[->] (A) to node[left] {$\phi_\alpha$} (C);
\draw[->] (B) to node[right] {$F$} (D);
\draw[->] (C) to node[below] {$F$} (D);

\node at (2,-1.5) {$\Downarrow \psi$};

\end{tikzpicture}
\]
of the following form.

We write $(D,F,\psi)$ to denote a $(y,\alpha)$-coequalizing square of the form above.
\end{definition}

\begin{definition}
Let $(D,F,\psi)$ and $(D',F',\psi')$ be $(y,\alpha)$-coequalizing squares. A \emph{morphism}
\[
H : (D,F,\psi) \to (D',F',\psi')
\]
consists of a functor $H : D \to D'$ such that
\[
H \circ F = F'
\]
and
\[
H \psi = \psi'.
\]
\end{definition}

We write $\Sq_{y,\alpha}$ for the category of $(y,\alpha)$-extension squares.

\begin{lemma}
The category $\Sq_{y,\alpha}$ has an initial object
\[
(E[y;\alpha], i_{y,\alpha}, \gamma).
\]
\end{lemma}

\begin{proof}
The $2$-category $\mathbf{Cat}$ has all coinserters (see Example 6.5 of \cite{FlexLim}). The initial object of 
$\Sq_{y,\alpha}$ is a specific example of a coinserter.
\end{proof}

\begin{lemma}\label{uni_prop_of_coins_clean}
For any $F : E \to D$, there is a natural bijection
\[
E/\mathbf{Cat}((E[y;\alpha], i_{y,\alpha}), (D,F))
\cong
[I_\alpha, D](F \circ *_y, F \circ \phi_\alpha).
\]
\end{lemma}

\begin{lemma}\label{coins_is_small_clean}
The morphism $i_{y,\alpha} : E \to E[y;\alpha]$ is a small object in $E/\mathbf{Cat}$.
\end{lemma}

\begin{proof}
The proof follows by expressing hom-sets via Lemma \ref{uni_prop_of_coins_clean}, using that $[I_\alpha,-]$ preserves filtered colimits, and implementing the fact that the forgetful functor 
\[
U:E/\mathbf{Cat} \to \mathbf{Cat}
\]
preserves filtered colimits.
\end{proof}

\begin{notation}
Let $E[y;\alpha]'$ denote the category obtained from $E[y;\alpha]$ by freely adjoining a morphism
\[
f_{y,\alpha} : y \to x_\alpha
\]
subject to the relation
\[
(\delta_\alpha)_i \circ f_{y,\alpha} = \gamma_i
\]
for every object $i$ of $I_\alpha$. The category $E[y;\alpha]'$ comes equipped with a functor 
\[
r_{y,\alpha}:E[y;\alpha]\to E[y;\alpha]'.
\]
\end{notation}

\begin{lemma}
Given a category $D$ together with a functor 
	\[
	F:E[y;\alpha]\to D
	\]
and a choice of morphism $h:F(y)\to F(x_\alpha)$ in $D$ such that 
	\[
	F((\delta_\alpha)_i)\circ h=F(\gamma_i)
	\]
	for every object $i$ of $I_\alpha$,
there is a unique functor 
	\[
	F':E[y;\alpha]'\to D
	\]
	such that $F'\circ r_{y,\alpha}=F$ and $F'(f_{y,\alpha})=h$.
\end{lemma}

\begin{notation}
We write $j_{y,\alpha}:E\to E[y;\alpha]'$ to be the following composite.
\[
j_{y,\alpha}:=r_{y,\alpha}\circ i_{y,\alpha}
\]
\end{notation}

\begin{lemma}\label{univ_prop_ind_map_is_fun_clean}
There is a bijection
\[
E/\mathbf{Cat}((E[y;\alpha]', j_{y,\alpha}), (D,F))
\cong
E/\mathbf{Cat}((E[y;\alpha], i_{y,\alpha}), (D,F))
\times
\mathbf{Cat}(F(y), F(x_\alpha))
\]
natural in the pair $(D,F)$.
\end{lemma}

\begin{lemma}\label{j_small_clean}
The functor $j_{y,\alpha} : E \to E[y;\alpha]'$ is a small object in $E/\mathbf{Cat}$.
\end{lemma}

\begin{notation}
Let $\Xi_E$ denote the set of morphisms in $E/\mathbf{Cat}$ that consists of maps of the following form in $E/\mathbf{Cat}$ where $y$ is an object of $ E$ and $\alpha\in A$.	
\[
\begin{tikzpicture}[node distance=2cm]
	\node (A){$E$};
	\node (B)[below of=A]{};
	\node (C)[left of=B]{$E[y;\alpha]$};
	\node (D)[node distance=4cm,right of=C]{$E[y;\alpha]'$};
	\draw[->](A) to node[left=3]{$i_{y,\alpha}$} (C);
	\draw[->](A) to node[right=3]{$j_{y,\alpha}$}(D);
	\draw[->](C) to node[below=3]{$r_{y,\alpha}$}(D);
\end{tikzpicture}
\]
Moreover, let
\[
(\mathcal{L}_E, \mathcal{R}_E), \quad (L_E, E_E, R_E)
\]
be the strong factorization system generated by $\Xi_E$ via Kelly’s Small Object Argument.
\end{notation}

\section{Universal Realization of Limit Sketches}

We now construct the universal realization functor for limit sketches and establish its universal property.

\begin{definition}\label{univ_lim_sk_const}
Define a $2$-functor
\[
\free : \limsk \to \relimsk
\]
on objects as follows. For a limit sketch $(E,A,\{C_\alpha\}_{\alpha \in A})$, set
\[
\free(E,A,\{C_\alpha\}) := (\mathrm{Fib}_E(!), A, \{C_\alpha'\}),
\]
where $\mathrm{Fib}_E(!)$ denotes the middle object in the factorization
\[
E \xrightarrow{L_E(!)} \mathrm{Fib}_E(!) \to *
\]
in $E/\mathbf{Cat}$, and where $C_\alpha'$ denotes the image of $C_\alpha$ under $L_E(!)$.
\end{definition}

\begin{lemma}
For each $\alpha \in A$, the cone $L_E(!)(C_\alpha)$ is a limit cone in $\mathrm{Fib}_E(!)$.
\end{lemma}

\begin{proof}
	Suppose $y\in\ob \mathrm{Fib}_E(!)$ and $\kappa:*_y\Rightarrow \phi_\alpha$ is a natural transformation. Then there is an induced functor 
	\[
	\Tilde{\kappa}:E[y;\alpha]\to \text{Fib}_E(!)
	\]
	defined by setting
\begin{itemize}
\item $\Tilde{\kappa}\circ i_{y,\alpha}=L_E(!)$ and
\item $\Tilde{\kappa}(\gamma_i)=\kappa_i$ for $i\in\ob(I_\alpha)$.
\end{itemize}
This induces the following commutative square.
	
	\[
	\begin{tikzpicture}[node distance=2.5cm,auto]
		\node (A) {$E[y;\alpha]$};
		\node (B)[right of=A] {$\mathrm{Fib}_E(!)$};
		\node(C)[below of=A] {$E[y;\alpha]'$};
		\node (D)[right of=C] {$*$};
		\draw[->] (A) to node[above=3] {$\Tilde{\kappa}$}(B);
		\draw[->,swap] (A) to node[left=3]{$r_{y,\alpha}$}(C);
		\draw[->] (B) to node[right=3]{$!$}(D);
		\draw[->,swap] (C) to node[below=3] {$!$}(D);
	\end{tikzpicture}.
	\]
	 By Kelly's Small Object Argument, there is a unique lift 
	\[
	w:E[y;\alpha]'\to \text{Fib}_E(!).
	\]
	 This means there is a uniquely induced pullback map from $(y,I_\alpha,\phi_\alpha,\kappa)$ to $C_\alpha$. Therefore $C_\alpha$ is sent to a limit cone in $\mathrm{Fib}_E(!)$.
\end{proof}

\begin{lemma}\label{uni_prop_of_real}
Let $F : E \to D$ be a functor sending each $C_\alpha$ to a limit cone in $D$. Then there exists a unique extension
\[
F' : \mathrm{Fib}_E(!) \to D
\]
such that
\[
F' \circ L_E(!) = F.
\]
\end{lemma}

\begin{proof}
Since $F:E\to D$ sends the cone $C_\alpha$ to a limit cone in $D$ for all $\alpha\in A$, the map $F:E\to D$ is a fibrant object of $E/\mathbf{Cat}$. This means that the terminal map
\[
\begin{tikzpicture}[node distance=2cm]
	\node (A){$E$};
	\node (B)[below of=A]{};
	\node (C)[left of=B]{$D$};
	\node (D)[node distance=4cm,right of=C]{$*$};
	\draw[->](A) to node[left=3]{$F$} (C);
	\draw[->](A) to node[right=3]{$!$}(D);
	\draw[->](C) to node[below=3]{$!$}(D);
\end{tikzpicture}
\]
is an element of $\mathcal{R}_E$. As $L_E(!):E\to \mathrm{Fib}_E(!)$ is an element of $\mathcal{R}_E$ and the diagram
\[
\begin{tikzpicture}[node distance=2.5cm,auto]
		\node (A) {$E$};
		\node (B)[right of=A] {$D$};
		\node(C)[below of=A] {$\mathrm{Fib}_E(!)$};
		\node (D)[right of=C] {$*$};
		\draw[->] (A) to node[above=3] {$F$}(B);
		\draw[->,swap] (A) to node[left=3]{$L_E(!)$}(C);
		\draw[->] (B) to node[right=3]{$!$}(D);
		\draw[->,swap] (C) to node[below=3] {$!$}(D);
	\end{tikzpicture}
\]
commutes, there is a unique functor $F':\mathrm{Fib}_E(!)\to D$ such that
	\[
	F'\circ L_E(!)=F.
	\]
    Therefore the map $L_E(!):E\to \mathrm{Fib}_E(!)$ sends the cone $C_\alpha$
	to a limit cone in $\mathrm{Fib}_E(!)$ for all $\alpha\in A.$
\end{proof}

\begin{lemma}\label{uni_prop_of_real_nat}
Let $F,G : E \to D$ be functors that preserve the specified limit cones and let $\phi : F \Rightarrow G$ be a natural transformation. Then there exists a unique natural transformation
\[
\phi' : F' \Rightarrow G'
\]
such that $\phi' L_E(!) = \phi$.
\end{lemma}

\begin{proof}
This follows by interpreting natural transformations as objects of $[P_0,[E,D]]$ and applying Lemma \ref{uni_prop_of_real} pointwise.
\end{proof}

\begin{lemma}
The assignment $\free$ extends to a $2$-functor
\[
\free : \limsk \to \relimsk
\]
and the maps $L_E(!)$ assemble into a $2$-natural transformation
\[
\eta : 1_{\limsk} \Rightarrow U \free.
\]
\end{lemma}

\begin{proof}
Functoriality on 1- and 2-cells follows from the universal property in Lemma \ref{uni_prop_of_real} and Lemma \ref{uni_prop_of_real_nat}.
\end{proof}

\begin{lemma}
If $E$ is already a realized limit sketch, then the canonical map
\[
L_E(!) : E \to \free(E)
\]
is an isomorphism.
\end{lemma}

\begin{proof}
Since all cones are already limits, $E$ satisfies the same universal property characterizing $\free(E)$.
\end{proof}

\begin{theorem}
The $2$-functor
\[
\free : \limsk \to \relimsk
\]
is left $2$-adjoint to the forgetful functor
\[
U : \relimsk \to \limsk.
\]
\end{theorem}

\begin{proof}
For any realized sketch $R$ and sketch $E$, morphisms
\[
\free(E) \to R
\]
correspond uniquely to morphisms
\[
E \to U(R)
\]
by Lemma \ref{uni_prop_of_real}. Naturality in both variables gives the required $2$-adjunction structure.
\end{proof}

\begin{remark}
This construction parallels normalization procedures for left sketches of \cite{Lob1}, but differs in that normalized left sketches are cocomplete. However, there is a $2$-functor from normalized left sketches to realized limit sketches where we forget about all the colimit cones. This allows us to view normalized left sketches as a reflective sub-$2$-category of realized limit sketches.
\end{remark}

\section{Models over Limit Sketches}

We conclude by interpreting the $2$-adjunction in terms of models and Morita equivalence.

\begin{definition}
Let $D$ be a complete category. A \emph{model} of a limit sketch $(E,A,\{C_\alpha\})$ in $D$ is a functor
\[
F : E \to D
\]
sending each $C_\alpha$ to a limit cone in $D$.
\end{definition}

\begin{notation}
We write $\mathrm{Mod}(E;D)$ for the full subcategory of $[E,D]$ consisting of models.
\end{notation}

\begin{definition}
A morphism of limit sketches $(F,f)$ is a \emph{Morita equivalence} if for every complete category $D$, the induced functor
\[
(F,f)^* : \mathrm{Mod}(E';D) \to \mathrm{Mod}(E;D)
\]
is an equivalence of categories.
\end{definition}

\begin{lemma}
For every limit sketch $E$, the unit map
\[
L_E(!) : E \to \free(E)
\]
is a Morita equivalence.
\end{lemma}

\begin{proof}
By the universal property of $\free(E)$, models of $E$ in any complete category $D$ are in canonical bijection with models of $\free(E)$ in $D$. This identification rises to an equivalence of categories.
\end{proof}

\end{document}